\documentclass[a4paper,leqno,12pt]{article}

\usepackage{amsmath}
\usepackage{amssymb}
\usepackage{amsthm}

\begin{document}

\bibliographystyle{amsalpha}
\newcommand\lra{\longrightarrow}
\newcommand\ra{\rightarrow}
\newtheorem{Problem}{Problem}[section]
\newtheorem{Assumption}{Assumption}[section]
\newtheorem{Theorem}{Theorem}[section]
\newtheorem{Lemma}{Lemma}[section]
\newtheorem{Remark}{Remark}[section]
\newtheorem{Corollary}{Corollary}[section]
\newtheorem{Conjecture}{Conjecture}[section]
\newtheorem{Proposition}{Proposition}[section]
\newtheorem{Example}{Example}[section]
\newtheorem{Definition}{Definition}[section]
\newtheorem{Question}{Question}[section]
\renewcommand{\thesubsection}{\it}

\title {Four-dimensional Painlev\'e systems of types $D_5^{(1)}$ and $B_4^{(1)}$ \\}

\author{Yusuke Sasano}
\maketitle

\begin{abstract}
We find and study a five-parameter family of four-dimensional coupled Painlev\'e V systems with affine Weyl group symmetry of type $D_5^{(1)}$. We then give an explicit description of a confluence from those systems to a four-parameter family of four-dimensional coupled Painlev\'e III systems with affine Weyl group symmetry of type $B_4^{(1)}$. 
\end{abstract}

\section{Introduction}

It is well-known that the Painlev\'e systems $P_{II},P_{III},P_{IV},P_{V}$ and $P_{VI}$ admit the affine Weyl groups of type $A_1^{(1)},C_2^{(1)},A_2^{(1)},A_3^{(1)}$ and $D_4^{(1)}$, respectively, as groups of B{\"a}cklund transformations. This suggests the following general problem (see \cite{N3}):
\begin{Problem}
For each affine root system $A$ with affine Weyl group $W(A)$, find a system of differential equations for which $W(A)$ acts as its B{\"a}cklund transformations.
\end{Problem}

One could expect that such nonlinear differential systems with affine Weyl group symmetry should admit rich mathematical structures, comparable with those of the Painlev\'e equations. In the case of type $A_l^{(1)}$, such equations are proposed in \cite{N1}. They are considered to be higher order versions of $P_V$ (resp. $P_{IV}$) when $l$ is odd (resp. even). These two examples by Noumi and Yamada motivated the author to find examples of higher order versions other than the systems of type $A_l^{(1)}$.

We will complete the study of the above problem in a series of four papers, for which this paper is the second, resulting in a series of equations for the remaining affine root systems of types $B_l^{(1)},C_l^{(1)}$ and $D_l^{(1)}$ (see \cite{Sasa4,Sasa5,Sasa6}). This paper is the stage in this project where we find and study four-dimensional coupled Painlev\'e V systems with $W(D_5^{(1)})$-symmetry explicitly given by
\begin{equation*}
\frac{dx}{dt}=\frac{\partial H_{D_5^{(1)}}}{\partial y}, \ \ \frac{dy}{dt}=-\frac{\partial H_{D_5^{(1)}}}{\partial x}, \ \ \frac{dz}{dt}=\frac{\partial H_{D_5^{(1)}}}{\partial w}, \ \ \frac{dw}{dt}=-\frac{\partial H_{D_5^{(1)}}}{\partial z}
\end{equation*}
with the Hamiltonian
\begin{align}\label{HD5}
\begin{split}
H_{D_5^{(1)}}=&H_{V}(x,y,t;\alpha_2+\alpha_5,\alpha_1,\alpha_2+2\alpha_3+\alpha_4)+H_{V}(z,w,t;\alpha_5,\alpha_3,\alpha_4)\\
&+\frac{2yz\{(z-1)w+\alpha_3\}}{t},
\end{split}
\end{align}
where the symbol $H_{V}(q,p,t;\gamma_1,\gamma_2,\gamma_3)$ denotes the Hamiltonian of the second-order Painlev\'e V systems given by
$$
H_{V}(q,p,t;\gamma_1,\gamma_2,\gamma_3)=\frac{q(q-1)p(p+t)-(\gamma_1+\gamma_3)qp+\gamma_1p+\gamma_2tq}{t}.
$$
Here $x,y,z$ and $w$ denote unknown complex variables, and $\alpha_0,\alpha_1,\dots,\alpha_5$ are complex parameters satisfying the relation:
$$
\alpha_0+\alpha_1+2\alpha_2+2\alpha_3+\alpha_4+\alpha_5=1.
$$
This is the first example which gave higher-order Painlev\'e equations of type $D_{2l+3}^{(1)}$. We then give an explicit description of a confluence from those systems to four-dimensional coupled Painlev\'e III systems with $W(B_4^{(1)})$-symmetry explicitly given by
\begin{equation*}
\frac{dx}{dt}=\frac{\partial H_{B_4^{(1)}}}{\partial y}, \ \ \frac{dy}{dt}=-\frac{\partial H_{B_4^{(1)}}}{\partial x}, \ \ \frac{dz}{dt}=\frac{\partial H_{B_4^{(1)}}}{\partial w}, \ \ \frac{dw}{dt}=-\frac{\partial H_{B_4^{(1)}}}{\partial z}
\end{equation*}
with the Hamiltonian
\begin{align}\label{HB4}
\begin{split}
H_{B_4^{(1)}} &=H_{III}(x,y,t;\alpha_0,\alpha_1)+H_{III}(z,w,t;\alpha_0+\alpha_1+2\alpha_2+\alpha_3,\alpha_3)\\
&+\frac{2yz(zw+\alpha_3)}{t},
\end{split}
\end{align}
where the symbol $H_{III}(q,p,t;\gamma_0,\gamma_1,\gamma_2)$ denotes the Hamiltonian of the second-order Painlev\'e III systems given by
$$
H_{III}(q,p,t;\gamma_0,\gamma_1,\gamma_2)=\frac{q^2p(p-1)+q\{(\gamma_0+\gamma_1)p-\gamma_1\}+tp}{t} \quad (\gamma_0+2\gamma_1+\gamma_2=1).
$$
Here $x,y,z$ and $w$ denote unknown complex variables and $\alpha_0,\alpha_1,\dots,\alpha_4$ are complex parameters satisfying the relation:
$$
\alpha_0+\alpha_1+2\alpha_2+2\alpha_3+2\alpha_4=1.
$$
This is the first example which gave higher-order Painlev\'e equations of type $B_{2l+2}^{(1)}$. We note that in \cite{Sasa5} we presented four-dimensional coupled Painlev\'e VI systems with $W(D_6^{(1)})$-symmetry. Before giving proofs (from Section 4 onward), we first state our results in the next three introductory sections. In Sections 1 and 2, we will present four-dimensional coupled Painlev\'e V and III systems with $W(D_5^{(1)})$-symmetry and $W(B_4^{(1)})$-symmetry, respectively. In Section 3, we will present four-dimensional polynomial Hamiltonian systems with $W(D_4^{(2)})$-symmetry. This is the first example which gave higher-order Painlev\'e equations of type $X_{l}^{(2)}$. In Section 4, we will consider the degeneration process from the system of type $D_5^{(1)}$ to the system of types $A_4^{(1)}$ and $B_4^{(1)}$, respectively. We also show that the B{\"a}cklund transformation groups for each root system are obtained from that for type $D_5^{(1)}$ by each degeneration process. After we review the notion of accessible singularity in Section 5, in Sections 6 and 7, we will make canonical coordinate systems for the system \eqref{eq:1} and \eqref{eq:2}, respectively.

\section{Main results for the case of $D_5^{(1)}$}

In this paper, we present a 5-parameter family of polynomial Hamiltonian systems that can be considered as four-dimensional coupled Painlev\'e V systems explicitly given by
\begin{equation}\label{eq:1}
  \left\{
  \begin{aligned}
   \frac{dx}{dt} =&\frac{2x^2y}{t}+x^2-\frac{2xy}{t}-\left(1+\frac{2\alpha_2+2\alpha_3+\alpha_5+\alpha_4}{t} \right)x\\
   &+\frac{\alpha_2+\alpha_5}{t}+\frac{2z((z-1)w+\alpha_3)}{t},\\
   \frac{dy}{dt} =&-\frac{2xy^2}{t}+\frac{y^2}{t}-2xy+\left(1+\frac{2\alpha_2+2\alpha_3+\alpha_5+\alpha_4}{t} \right)y-\alpha_1,\\
   \frac{dz}{dt} =&\frac{2z^2w}{t}+z^2-\frac{2zw}{t}-\left(1+\frac{\alpha_5+\alpha_4}{t} \right)z+\frac{\alpha_5}{t}+\frac{2yz(z-1)}{t},\\
   \frac{dw}{dt} =&-\frac{2zw^2}{t}+\frac{w^2}{t}-2zw+\left(1+\frac{\alpha_5+\alpha_4}{t} \right)w-\alpha_3-\frac{2y(-w+2zw+\alpha_3)}{t}
   \end{aligned}
  \right. 
\end{equation}
with the Hamiltonian \eqref{HD5}.

\begin{Proposition}\label{prop:1}
The system \eqref{eq:1} has the following invariant divisors\rm{:\rm}
\begin{center}
\begin{tabular}{|c|c|c|} \hline
codimension & invariant divisors & parameter's relation \\ \hline
1 & $f_0:=y+t$ & $\alpha_0=0$  \\ \hline
1 & $f_1:=y$ & $\alpha_1=0$  \\ \hline
1 & $f_2:=x-z$ & $\alpha_2=0$  \\ \hline
1 & $f_3:=w$ & $\alpha_3=0$  \\ \hline
1 & $f_4:=z-1$ & $\alpha_4=0$   \\ \hline
1 & $f_5:=z$ & $\alpha_5=0$   \\ \hline
\end{tabular}
\end{center}
\end{Proposition}

\noindent
The list must be read as follows. Setting $\alpha_1=0$, then the system \eqref{eq:1} admits a particular solution $y=0$. Moreover $(z,w)$ satisfy the fifth Painlev\'e system. And $x$ satisfies Riccati equations whose coefficients are polynomials in $(z,w)$, and so on.

\begin{figure}[ht]
\unitlength 0.1in
\begin{picture}(49.80,20.80)(8.10,-23.70)
%
\special{pn 20}%
\special{ar 2170 842 240 195  1.5707963 6.2831853}%
\special{ar 2170 842 240 195  0.0000000 1.5291537}%
%
\special{pn 20}%
\special{ar 2170 1766 240 194  1.5707963 6.2831853}%
\special{ar 2170 1766 240 194  0.0000000 1.5291537}%
%
\special{pn 20}%
\special{ar 2960 1320 240 195  1.5707963 6.2831853}%
\special{ar 2960 1320 240 195  0.0000000 1.5291537}%
%
\special{pn 20}%
\special{ar 4790 845 240 195  1.5707963 6.2831853}%
\special{ar 4790 845 240 195  0.0000000 1.5291537}%
%
\special{pn 20}%
\special{ar 4790 1777 240 194  1.5707963 6.2831853}%
\special{ar 4790 1777 240 194  0.0000000 1.5291537}%
%
\special{pn 20}%
\special{pa 2380 924}%
\special{pa 2770 1183}%
\special{fp}%
%
\special{pn 20}%
\special{pa 2390 1693}%
\special{pa 2790 1458}%
\special{fp}%
%
\special{pn 20}%
\special{ar 3880 1328 240 195  1.5707963 6.2831853}%
\special{ar 3880 1328 240 195  0.0000000 1.5291537}%
%
\special{pn 20}%
\special{pa 3200 1328}%
\special{pa 3620 1328}%
\special{fp}%
%
\special{pn 20}%
\special{pa 4110 1210}%
\special{pa 4580 967}%
\special{fp}%
%
\special{pn 20}%
\special{pa 4080 1453}%
\special{pa 4560 1712}%
\special{fp}%
%
\put(15.6000,-8.9900){\makebox(0,0)[lb]{}}%
\put(46.1000,-9.5000){\makebox(0,0)[lb]{$z-1$}}%
\put(37.4000,-14.3000){\makebox(0,0)[lb]{$w$}}%
\put(27.5000,-13.9000){\makebox(0,0)[lb]{$x-z$}}%
\put(19.9000,-18.4000){\makebox(0,0)[lb]{$y+t$}}%
\put(20.6000,-9.3000){\makebox(0,0)[lb]{$y$}}%
\put(24.4000,-20.4000){\makebox(0,0)[lb]{Dynkin diagram of type ${D_5}^{(1)}$}}%
\put(46.7000,-18.7000){\makebox(0,0)[lb]{$z$}}%
%
\special{pn 8}%
\special{pa 810 1340}%
\special{pa 5790 1330}%
\special{dt 0.045}%
\special{pa 5790 1330}%
\special{pa 5789 1330}%
\special{dt 0.045}%
%
\special{pn 8}%
\special{pa 3410 290}%
\special{pa 3410 2370}%
\special{dt 0.045}%
\special{pa 3410 2370}%
\special{pa 3410 2369}%
\special{dt 0.045}%
%
\special{pn 20}%
\special{pa 1170 1290}%
\special{pa 1138 1266}%
\special{pa 1107 1247}%
\special{pa 1079 1237}%
\special{pa 1055 1243}%
\special{pa 1037 1265}%
\special{pa 1025 1299}%
\special{pa 1019 1340}%
\special{pa 1018 1383}%
\special{pa 1023 1423}%
\special{pa 1034 1453}%
\special{pa 1051 1470}%
\special{pa 1073 1470}%
\special{pa 1099 1455}%
\special{pa 1128 1431}%
\special{pa 1159 1401}%
\special{pa 1160 1400}%
\special{sp}%
%
\special{pn 20}%
\special{pa 3360 590}%
\special{pa 3325 574}%
\special{pa 3292 557}%
\special{pa 3266 538}%
\special{pa 3247 517}%
\special{pa 3239 492}%
\special{pa 3243 464}%
\special{pa 3257 435}%
\special{pa 3281 408}%
\special{pa 3310 384}%
\special{pa 3345 367}%
\special{pa 3383 358}%
\special{pa 3422 359}%
\special{pa 3460 368}%
\special{pa 3494 384}%
\special{pa 3522 405}%
\special{pa 3541 430}%
\special{pa 3550 457}%
\special{pa 3546 486}%
\special{pa 3532 515}%
\special{pa 3513 545}%
\special{pa 3510 550}%
\special{sp}%
%
\special{pn 20}%
\special{pa 1100 1440}%
\special{pa 1200 1360}%
\special{fp}%
\special{sh 1}%
\special{pa 1200 1360}%
\special{pa 1135 1386}%
\special{pa 1158 1393}%
\special{pa 1160 1417}%
\special{pa 1200 1360}%
\special{fp}%
%
\special{pn 20}%
\special{pa 3550 460}%
\special{pa 3460 640}%
\special{fp}%
\special{sh 1}%
\special{pa 3460 640}%
\special{pa 3508 589}%
\special{pa 3484 592}%
\special{pa 3472 571}%
\special{pa 3460 640}%
\special{fp}%
\special{pa 3460 640}%
\special{pa 3460 640}%
\special{fp}%
\put(8.5000,-12.2000){\makebox(0,0)[lb]{${\pi}_1$}}%
\put(35.9000,-6.0000){\makebox(0,0)[lb]{${\pi}_2$}}%
\put(18.1000,-15.8000){\makebox(0,0)[lb]{$\alpha_0$}}%
\put(18.1000,-6.7000){\makebox(0,0)[lb]{$\alpha_1$}}%
\put(27.5000,-11.1000){\makebox(0,0)[lb]{$\alpha_2$}}%
\put(37.0000,-11.0000){\makebox(0,0)[lb]{$\alpha_3$}}%
\put(46.0000,-6.3000){\makebox(0,0)[lb]{$\alpha_4$}}%
\put(46.3000,-15.7000){\makebox(0,0)[lb]{$\alpha_5$}}%
\end{picture}%
\label{CPV1}
\caption{The transformations described in Theorem \ref{th:1.1} define a representation of the affine Weyl group of type $D_5^{(1)}$, that is, they satisfy the following relations: ${s_0}^2={s_1}^2=\dots={s_5}^2={\pi_1}^2=\dots={\pi_4}^2=1, \ (s_0s_1)^2=(s_0s_3)^2=(s_0s_4)^2=(s_0s_5)^2=(s_1s_3)^2=(s_1s_4)^2=(s_1s_5)^2=(s_2s_4)^2=(s_2s_5)^2=(s_4s_5)^2=1, \ (s_0s_2)^3=(s_1s_2)^3=(s_2s_3)^3=(s_3s_4)^3=(s_3s_5)^3=1$. The symbol in each circle denotes the invariant divisor $f_i$ of the system \eqref{eq:1} (see Proposition \ref{prop:1}).}
\end{figure}
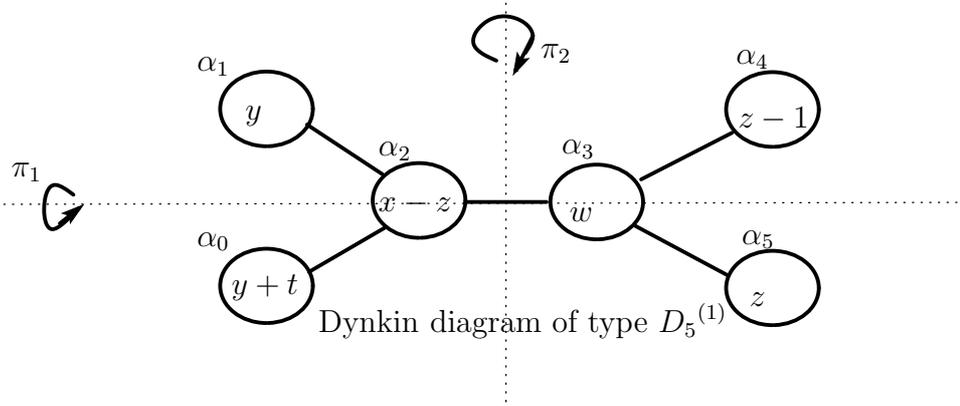

\begin{Theorem}\label{th:1.1}
The system \eqref{eq:1} admits extended affine Weyl group symmetry of type $D_5^{(1)}$ as the group of its B{\"a}cklund transformations, whose generators $s_0,s_1,\dots,s_5,{\pi}_1,\\
{\pi}_2,\dots,{\pi}_4$ defined as follows$:$ with {\it the notation} $(*):=(x,y,z,w,t;\alpha_0,\alpha_1,\dots,\alpha_5)$,
\begin{align}
\begin{split}
s_0: (*) \rightarrow &\left(x+\frac{\alpha_0}{y+t},y,z,w,t;-\alpha_0,\alpha_1,\alpha_2+\alpha_0,\alpha_3,\alpha_4,\alpha_5 \right),\\
s_1: (*) \rightarrow &\left(x+\frac{\alpha_1}{y},y,z,w,t;\alpha_0,-\alpha_1,\alpha_2+\alpha_1,\alpha_3,\alpha_4,\alpha_5 \right),\\
s_2: (*) \rightarrow &\left(x,y-\frac{\alpha_2}{x-z},z,w+\frac{\alpha_2}{x-z},t;\alpha_0+\alpha_2,\alpha_1+\alpha_2,-\alpha_2,\alpha_3+\alpha_2,\alpha_4,\alpha_5 \right),\\
s_3: (*) \rightarrow &\left(x,y,z+\frac{\alpha_3}{w},w,t;\alpha_0,\alpha_1,\alpha_2+\alpha_3,-\alpha_3,\alpha_4+\alpha_3,\alpha_5+\alpha_3 \right),\\
s_4: (*) \rightarrow &\left(x,y,z,w-\frac{\alpha_4}{(z-1)},t;\alpha_0,\alpha_1,\alpha_2,\alpha_3+\alpha_4,-\alpha_4,\alpha_5 \right),\\
s_5: (*) \rightarrow &\left(x,y,z,w-\frac{\alpha_5}{z},t;\alpha_0,\alpha_1,\alpha_2,\alpha_3+\alpha_5,\alpha_4,-\alpha_5 \right),\\
\pi_{1}: (*) \rightarrow &(1-x,-y-t,1-z,-w,t;\alpha_1,\alpha_0,\alpha_2,\alpha_3,\alpha_5,\alpha_4),\\
\pi_{2}: (*) \rightarrow &\left((y+w+t)/t,-t(z-1),(y+t)/t,-t(x-z),-t;\alpha_5,\alpha_4,\alpha_3,\alpha_2,\alpha_1,\alpha_0 \right),\\
\pi_{3}: (*) \rightarrow &(1-x,-y,1-z,-w,-t;\alpha_0,\alpha_1,\alpha_2,\alpha_3,\alpha_5,\alpha_4),\\
\pi_{4}: (*) \rightarrow &(x,y+t,z,w,-t;\alpha_1,\alpha_0,\alpha_2,\alpha_3,\alpha_4,\alpha_5).
\end{split}
\end{align}
\end{Theorem}

We note that it is easy to see that the generators ${\pi}_2,{\pi}_3,{\pi}_4$ satisfy the relation$:$
$$
{\pi}_4={\pi}_2{\pi}_3{\pi}_2.
$$

We remark that the B{\"a}cklund transformations of the system \eqref{eq:1} have the universal description to root system of type $D_5^{(1)}$ (see \cite{N2}). Since this universal B{\"a}cklund transformations have Lie theoritic origin, similarity reduction of Drinfeld-Sokolov hierarchy admits such B{\"a}cklund symmetry.

\begin{Remark}
Y. Yamada found an autonomous version of this system.
\end{Remark}

\begin{Proposition}
Let us define the following translation operators
$$
T_1:={{\pi}_1}s_5s_3s_2s_1s_0s_2s_3s_5, \ T_2:={{\pi}_2}T_1{{\pi}_2}, \ T_3:=s_1s_4T_1s_4s_1,
$$
$$
T_4:=s_2s_3T_3s_3s_2, \ T_5:=s_1T_4s_1, \ T_6:=s_3T_3s_3.
$$
These translation operators act on parameters $\alpha_i$ as follows$:$
\begin{align}
\begin{split}
T_1(\alpha_0,\alpha_1,\dots,\alpha_5)=&(\alpha_0,\alpha_1,\dots,\alpha_5)+(0,0,0,0,1,-1),\\
T_2(\alpha_0,\alpha_1,\dots,\alpha_5)=&(\alpha_0,\alpha_1,\dots,\alpha_5)+(-1,1,0,0,0,0),\\
T_3(\alpha_0,\alpha_1,\dots,\alpha_5)=&(\alpha_0,\alpha_1,\dots,\alpha_5)+(0,0,0,1,-1,-1),\\
T_4(\alpha_0,\alpha_1,\dots,\alpha_5)=&(\alpha_0,\alpha_1,\dots,\alpha_5)+(1,1,-1,0,0,0),\\
T_5(\alpha_0,\alpha_1,\dots,\alpha_5)=&(\alpha_0,\alpha_1,\dots,\alpha_5)+(0,0,1,-1,0,0).
\end{split}
\end{align}
\end{Proposition}

\begin{figure}[h]
\unitlength 0.1in
\begin{picture}(53.70,27.46)(2.70,-29.07)
%
\special{pn 20}%
\special{ar 2160 1995 240 158  1.5707963 6.2831853}%
\special{ar 2160 1995 240 158  0.0000000 1.5291537}%
%
\special{pn 20}%
\special{ar 2160 2742 240 157  1.5707963 6.2831853}%
\special{ar 2160 2742 240 157  0.0000000 1.5291537}%
%
\special{pn 20}%
\special{ar 2950 2382 240 157  1.5707963 6.2831853}%
\special{ar 2950 2382 240 157  0.0000000 1.5291537}%
%
\special{pn 20}%
\special{ar 4780 1998 240 158  1.5707963 6.2831853}%
\special{ar 4780 1998 240 158  0.0000000 1.5291537}%
%
\special{pn 20}%
\special{ar 4780 2751 240 156  1.5707963 6.2831853}%
\special{ar 4780 2751 240 156  0.0000000 1.5291537}%
%
\special{pn 20}%
\special{pa 2370 2062}%
\special{pa 2760 2271}%
\special{fp}%
%
\special{pn 20}%
\special{pa 2380 2683}%
\special{pa 2780 2493}%
\special{fp}%
%
\special{pn 20}%
\special{ar 3870 2388 240 157  1.5707963 6.2831853}%
\special{ar 3870 2388 240 157  0.0000000 1.5291537}%
%
\special{pn 20}%
\special{pa 3190 2388}%
\special{pa 3610 2388}%
\special{fp}%
%
\special{pn 20}%
\special{pa 4100 2293}%
\special{pa 4570 2096}%
\special{fp}%
%
\special{pn 20}%
\special{pa 4070 2489}%
\special{pa 4550 2698}%
\special{fp}%
%
\put(49.8000,-11.8700){\makebox(0,0)[lb]{}}%
\put(45.8000,-21.0000){\makebox(0,0)[lb]{$z-1$}}%
\put(37.8000,-24.5900){\makebox(0,0)[lb]{$w$}}%
\put(27.7000,-24.5900){\makebox(0,0)[lb]{$x-z$}}%
\put(19.9000,-28.3000){\makebox(0,0)[lb]{$y+t$}}%
\put(20.5000,-21.0800){\makebox(0,0)[lb]{$y$}}%
\put(24.2000,-30.3800){\makebox(0,0)[lb]{Dynkin diagram of type ${D_5}^{(1)}$}}%
\put(46.7000,-28.3000){\makebox(0,0)[lb]{$z$}}%
%
\special{pn 20}%
\special{ar 4610 365 240 158  1.5707963 6.2831853}%
\special{ar 4610 365 240 158  0.0000000 1.5291537}%
%
\special{pn 20}%
\special{ar 4610 1111 240 157  1.5707963 6.2831853}%
\special{ar 4610 1111 240 157  0.0000000 1.5291537}%
%
\special{pn 8}%
\special{ar 5400 751 240 158  1.5707963 1.6310978}%
\special{ar 5400 751 240 158  1.8120024 1.8723039}%
\special{ar 5400 751 240 158  2.0532084 2.1135099}%
\special{ar 5400 751 240 158  2.2944144 2.3547159}%
\special{ar 5400 751 240 158  2.5356204 2.5959220}%
\special{ar 5400 751 240 158  2.7768265 2.8371280}%
\special{ar 5400 751 240 158  3.0180325 3.0783340}%
\special{ar 5400 751 240 158  3.2592385 3.3195400}%
\special{ar 5400 751 240 158  3.5004446 3.5607461}%
\special{ar 5400 751 240 158  3.7416506 3.8019521}%
\special{ar 5400 751 240 158  3.9828566 4.0431581}%
\special{ar 5400 751 240 158  4.2240627 4.2843642}%
\special{ar 5400 751 240 158  4.4652687 4.5255702}%
\special{ar 5400 751 240 158  4.7064747 4.7667762}%
\special{ar 5400 751 240 158  4.9476807 5.0079823}%
\special{ar 5400 751 240 158  5.1888868 5.2491883}%
\special{ar 5400 751 240 158  5.4300928 5.4903943}%
\special{ar 5400 751 240 158  5.6712988 5.7316003}%
\special{ar 5400 751 240 158  5.9125049 5.9728064}%
\special{ar 5400 751 240 158  6.1537109 6.2140124}%
\special{ar 5400 751 240 158  6.3949169 6.4552184}%
\special{ar 5400 751 240 158  6.6361230 6.6964245}%
\special{ar 5400 751 240 158  6.8773290 6.9376305}%
\special{ar 5400 751 240 158  7.1185350 7.1788365}%
\special{ar 5400 751 240 158  7.3597411 7.4200426}%
\special{ar 5400 751 240 158  7.6009471 7.6612486}%
%
\special{pn 8}%
\special{pa 4820 431}%
\special{pa 5210 641}%
\special{dt 0.045}%
\special{pa 5210 641}%
\special{pa 5209 641}%
\special{dt 0.045}%
%
\special{pn 8}%
\special{pa 4830 1052}%
\special{pa 5230 862}%
\special{dt 0.045}%
\special{pa 5230 862}%
\special{pa 5229 862}%
\special{dt 0.045}%
\put(44.3000,-4.3000){\makebox(0,0)[lb]{$z-1$}}%
\put(35.7000,-8.2500){\makebox(0,0)[lb]{$w$}}%
\put(45.1000,-12.0000){\makebox(0,0)[lb]{$z$}}%
%
\special{pn 20}%
\special{ar 2180 318 240 157  1.5707963 6.2831853}%
\special{ar 2180 318 240 157  0.0000000 1.5291537}%
%
\special{pn 20}%
\special{ar 2180 1070 240 157  1.5707963 6.2831853}%
\special{ar 2180 1070 240 157  0.0000000 1.5291537}%
%
\special{pn 8}%
\special{ar 1270 707 240 158  1.5707963 1.6310978}%
\special{ar 1270 707 240 158  1.8120024 1.8723039}%
\special{ar 1270 707 240 158  2.0532084 2.1135099}%
\special{ar 1270 707 240 158  2.2944144 2.3547159}%
\special{ar 1270 707 240 158  2.5356204 2.5959220}%
\special{ar 1270 707 240 158  2.7768265 2.8371280}%
\special{ar 1270 707 240 158  3.0180325 3.0783340}%
\special{ar 1270 707 240 158  3.2592385 3.3195400}%
\special{ar 1270 707 240 158  3.5004446 3.5607461}%
\special{ar 1270 707 240 158  3.7416506 3.8019521}%
\special{ar 1270 707 240 158  3.9828566 4.0431581}%
\special{ar 1270 707 240 158  4.2240627 4.2843642}%
\special{ar 1270 707 240 158  4.4652687 4.5255702}%
\special{ar 1270 707 240 158  4.7064747 4.7667762}%
\special{ar 1270 707 240 158  4.9476807 5.0079823}%
\special{ar 1270 707 240 158  5.1888868 5.2491883}%
\special{ar 1270 707 240 158  5.4300928 5.4903943}%
\special{ar 1270 707 240 158  5.6712988 5.7316003}%
\special{ar 1270 707 240 158  5.9125049 5.9728064}%
\special{ar 1270 707 240 158  6.1537109 6.2140124}%
\special{ar 1270 707 240 158  6.3949169 6.4552184}%
\special{ar 1270 707 240 158  6.6361230 6.6964245}%
\special{ar 1270 707 240 158  6.8773290 6.9376305}%
\special{ar 1270 707 240 158  7.1185350 7.1788365}%
\special{ar 1270 707 240 158  7.3597411 7.4200426}%
\special{ar 1270 707 240 158  7.6009471 7.6612486}%
%
\special{pn 8}%
\special{pa 1500 612}%
\special{pa 1970 416}%
\special{dt 0.045}%
\special{pa 1970 416}%
\special{pa 1969 416}%
\special{dt 0.045}%
%
\special{pn 8}%
\special{pa 1470 809}%
\special{pa 1950 1018}%
\special{dt 0.045}%
\special{pa 1950 1018}%
\special{pa 1949 1018}%
\special{dt 0.045}%
\put(28.7000,-8.2500){\makebox(0,0)[lb]{$x$}}%
\put(19.8000,-11.3800){\makebox(0,0)[lb]{$y+t$}}%
\put(20.5000,-4.1100){\makebox(0,0)[lb]{$y$}}%
%
\special{pn 20}%
\special{ar 3670 725 240 158  1.5707963 6.2831853}%
\special{ar 3670 725 240 158  0.0000000 1.5291537}%
%
\special{pn 20}%
\special{pa 3900 630}%
\special{pa 4370 434}%
\special{fp}%
%
\special{pn 20}%
\special{pa 3870 826}%
\special{pa 4350 1035}%
\special{fp}%
\put(52.1000,-8.5400){\makebox(0,0)[lb]{$w+t$}}%
%
\special{pn 20}%
\special{ar 2990 739 240 156  1.5707963 6.2831853}%
\special{ar 2990 739 240 156  0.0000000 1.5291537}%
%
\special{pn 20}%
\special{pa 2410 418}%
\special{pa 2800 627}%
\special{fp}%
%
\special{pn 20}%
\special{pa 2420 1040}%
\special{pa 2820 850}%
\special{fp}%
\put(10.9000,-8.0600){\makebox(0,0)[lb]{$x-1$}}%
%
\special{pn 20}%
\special{pa 3000 1346}%
\special{pa 3000 1652}%
\special{fp}%
\special{sh 1}%
\special{pa 3000 1652}%
\special{pa 3020 1585}%
\special{pa 3000 1599}%
\special{pa 2980 1585}%
\special{pa 3000 1652}%
\special{fp}%
%
\special{pn 20}%
\special{pa 3660 1346}%
\special{pa 3660 1652}%
\special{fp}%
\special{sh 1}%
\special{pa 3660 1652}%
\special{pa 3680 1585}%
\special{pa 3660 1599}%
\special{pa 3640 1585}%
\special{pa 3660 1652}%
\special{fp}%
\put(2.7000,-13.8000){\makebox(0,0)[lb]{Dynkin diagram of type ${A_3}^{(1)}$}}%
\end{picture}%
\label{CPV2}
\caption{We make the Dynkin diagram of type $D_5^{(1)}$ by connecting two Dynkin diagrams of type $A_3^{(1)}$ in addition to the term with invariant divisor $x-z$.}
\end{figure}
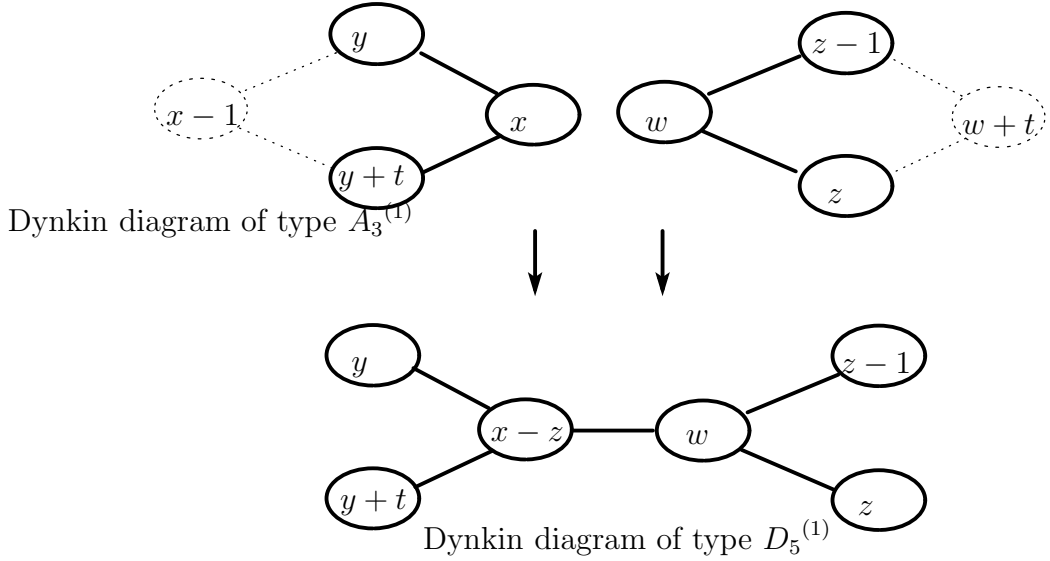

\begin{Theorem}\label{th:1.2}
Let us consider a polynomial Hamiltonian system with Hamiltonian $H \in {\Bbb C}(t)[x,y,z,w]$. We assume that

$(A1)$ $deg(H)=6$ with respect to $x,y,z,w$.

$(A2)$ This system becomes again a polynomial Hamiltonian system in each coordinate system $(x_i,y_i,z_i,w_i) \ (i=0,1,\dots ,5)${\rm : \rm}
\begin{align}
\begin{split}
&x_0=1/x, \ y_0=-((y+t)x+\alpha_0)x, \ z_0=z, \ w_0=w,\\
&x_1=1/x, \ y_1=-(yx+\alpha_1)x, \ z_1=z, \ w_1=w,\\
&x_2=-((x-z)y-\alpha_2)y, \ y_2=1/y, \ z_2=z, \ w_2=w+y,\\
&x_3=x, \ y_3=y, \ z_3=1/z, \ w_3=-(wz+\alpha_3)z, \\
&x_4=x, \ y_4=y, \ z_4=-((z-1)w-\alpha_4)w, \ w_4=1/w,\\
&x_5=x, \ y_5=y, \ z_5=-(zw-\alpha_5)w, \ w_5=1/w.
\end{split}
\end{align}
Then such a system coincides with the system \eqref{eq:1}.
\end{Theorem}

In addition to Theorems \ref{th:1.1} and \ref{th:1.2}, we will prove that the system \eqref{eq:1} degenerates to the system of type $A_4^{(1)}$ by taking a coupling confluence process of $P_{V} \ra P_{IV}$.
\begin{Theorem}\label{th:1.3}
For the system \eqref{eq:1} of type $D_5^{(1)}$, we make the change of parameters and variables$:$
\begin{gather}
\begin{gathered}\label{11}
\alpha_0=A_0-A_2-A_3+\frac{1}{2}{\varepsilon}^{-2}, \ \alpha_1=A_1, \ \alpha_2=A_2, \ \alpha_3=A_3, \ \alpha_4=-\frac{1}{2}{\varepsilon}^{-2}, \ \alpha_5=A_4,
\end{gathered}\\
\begin{gathered}\label{12}
t=\frac{1}{2}{\varepsilon}^{-2}(1+2\varepsilon T), \ x=-\frac{\varepsilon X}{1-\varepsilon X}, \ y=-\varepsilon^{-1}(1-\varepsilon X)[Y-\varepsilon (A_1+XY)],\\
z=-\frac{\varepsilon Z}{1-\varepsilon Z}, \ w=-\varepsilon^{-1}(1-\varepsilon Z)[W-\varepsilon (A_3+ZW)]
\end{gathered}
\end{gather}
from $\alpha_0,\alpha_1,\dots,\alpha_5,t,x,y,z,w$ to $A_0,A_1,\dots,A_4,\varepsilon,T,X,Y,Z,W$. Then the system \eqref{eq:1} can also be written in the new variables $T,X,Y,Z,W$ and parameters $A_0,A_1,\dots,\\
A_4,\varepsilon$ as a Hamiltonian system. This new system tends to the system of type $A_4^{(1)}$ as $\varepsilon \ra 0$.
\end{Theorem}

By proving the following theorem, we see how the degeneration process in Theorem \ref{th:1.3} works on the B{\"a}cklund transformation group $W(D_5^{(1)})=<s_0,s_1,\dots,s_5>$ described in Theorem \ref{th:1.1}.
\begin{Theorem}\label{D5A4}
For the degeneration process in Theorem \ref{th:1.3}, we can choose a subgroup $W_{D_5^{(1)} \ra A_4^{(1)}}$ of the B{\"a}cklund transformation group $W(D_5^{(1)})$ so that $W_{D_5^{(1)} \rightarrow A_4^{(1)}}$ converges to $W(A_4^{(1)})$ as $\varepsilon \rightarrow 0$.
\end{Theorem}

\section{Main results for the case of $B_4^{(1)}$}

It is well-known that the fifth Painlev\'e equation $P_{V}$ has a confluence to the third Painlev\'e equation $P_{III}$, where two accessible singularities come together into a single singularity. This suggests the possibility that there exists a procedure for searching for fourth-order versions of Painlev\'e III, by using Takano's description of the confluence process from $P_{V}$ to $P_{III}$ for the coordinate systems $(x,y)$ and $(z,w)$, respectively (see \cite{T2,T3}).  In this vein, the goal of this work is to find a fourth-order version of the Painlev\'e III equation with symmetry under the group which degenerates from the affine Weyl group of type $D_5^{(1)}$ by a coupling confluence process. In this paper, we also present a 4-parameter family of polynomial Hamiltonian systems that can be considered as four-dimensional coupled Painlev\'e III systems explicitly given by
\begin{equation}\label{eq:2}
  \left\{
  \begin{aligned}
   \frac{dx}{dt} &=\frac{2x^2y-{x^2}+(\alpha_0+\alpha_1)x+2\alpha_3z+2z^2w}{t}+1,\\
   \frac{dy}{dt} &=\frac{-2xy^2+2xy-(\alpha_0+\alpha_1)y+{\alpha_1}}{t},\\
   \frac{dz}{dt} &=\frac{2z^2w-{z^2}+(\alpha_0+\alpha_1+2\alpha_2+2\alpha_3)z+2yz^2}{t}+1,\\
   \frac{dw}{dt} &=\frac{-2zw^2+2{zw}-(\alpha_0+\alpha_1+2\alpha_2+2\alpha_3)w-2\alpha_3y-4yzw+\alpha_3}{t}
   \end{aligned}
  \right. 
\end{equation}
with the Hamiltonian \eqref{HB4}.

\begin{Proposition}\label{prop:2}
The system \eqref{eq:2} has the following invariant divisors\rm{:\rm}
\begin{center}
\begin{tabular}{|c|c|c|} \hline
codimension & invariant divisors & parameter's relation \\ \hline
1 & $f_0:=y-1$ & $\alpha_0=0$  \\ \hline
1 & $f_1:=y$ & $\alpha_1=0$  \\ \hline
1 & $f_2:=x-z$ & $\alpha_2=0$  \\ \hline
1 & $f_3:=w$ & $\alpha_3=0$  \\ \hline
\end{tabular}
\end{center}
\end{Proposition}

\begin{Theorem}\label{th:2.1}
The system \eqref{eq:2} admits extended affine Weyl group symmetry of type $B_4^{(1)}$ as the group of its B{\"a}cklund transformations, whose generators $s_0,s_1,\dots,s_4,{\pi}_1,\\
{\pi}_2$ defined as follows$:$ with {\it the notation} $(*):=(x,y,z,w,t;\alpha_0,\alpha_1,\dots,\alpha_4)$,

\begin{figure}[ht]
\unitlength 0.1in
\begin{picture}(44.89,17.03)(3.00,-21.03)
%
\put(3.0000,-9.8900){\makebox(0,0)[lb]{}}%
%
\special{pn 20}%
\special{ar 1710 820 189 223  1.5707963 6.2831853}%
\special{ar 1710 820 189 223  0.0000000 1.5284935}%
\put(13.8000,-5.7000){\makebox(0,0)[lb]{$\alpha_1$}}%
%
\special{pn 20}%
\special{ar 1710 1880 189 223  1.5707963 6.2831853}%
\special{ar 1710 1880 189 223  0.0000000 1.5284935}%
%
\special{pn 20}%
\special{ar 2680 1330 189 223  1.5707963 6.2831853}%
\special{ar 2680 1330 189 223  0.0000000 1.5284935}%
%
\special{pn 20}%
\special{ar 3630 1350 189 223  1.5707963 6.2831853}%
\special{ar 3630 1350 189 223  0.0000000 1.5284935}%
%
\special{pn 20}%
\special{ar 4600 1340 189 223  1.5707963 6.2831853}%
\special{ar 4600 1340 189 223  0.0000000 1.5284935}%
%
\special{pn 20}%
\special{pa 1880 920}%
\special{pa 2530 1230}%
\special{fp}%
%
\special{pn 20}%
\special{pa 1900 1860}%
\special{pa 2530 1480}%
\special{fp}%
%
\special{pn 20}%
\special{pa 2880 1330}%
\special{pa 3420 1330}%
\special{fp}%
%
\special{pn 20}%
\special{pa 3820 1280}%
\special{pa 4410 1280}%
\special{fp}%
\special{sh 1}%
\special{pa 4410 1280}%
\special{pa 4343 1260}%
\special{pa 4357 1280}%
\special{pa 4343 1300}%
\special{pa 4410 1280}%
\special{fp}%
%
\special{pn 20}%
\special{pa 3810 1410}%
\special{pa 4410 1410}%
\special{fp}%
\special{sh 1}%
\special{pa 4410 1410}%
\special{pa 4343 1390}%
\special{pa 4357 1410}%
\special{pa 4343 1430}%
\special{pa 4410 1410}%
\special{fp}%
\put(13.8000,-16.4000){\makebox(0,0)[lb]{$\alpha_0$}}%
\put(23.8000,-11.0000){\makebox(0,0)[lb]{$\alpha_2$}}%
\put(33.8000,-11.2000){\makebox(0,0)[lb]{$\alpha_3$}}%
\put(43.7000,-11.2000){\makebox(0,0)[lb]{$\alpha_4$}}%
%
\special{pn 20}%
\special{pa 1190 960}%
\special{pa 1162 978}%
\special{pa 1135 997}%
\special{pa 1108 1016}%
\special{pa 1082 1035}%
\special{pa 1057 1056}%
\special{pa 1033 1077}%
\special{pa 1010 1099}%
\special{pa 989 1123}%
\special{pa 970 1148}%
\special{pa 953 1175}%
\special{pa 938 1204}%
\special{pa 926 1234}%
\special{pa 915 1266}%
\special{pa 908 1300}%
\special{pa 902 1334}%
\special{pa 899 1369}%
\special{pa 898 1404}%
\special{pa 899 1440}%
\special{pa 902 1475}%
\special{pa 908 1510}%
\special{pa 916 1545}%
\special{pa 926 1578}%
\special{pa 937 1611}%
\special{pa 951 1642}%
\special{pa 968 1672}%
\special{pa 986 1700}%
\special{pa 1006 1726}%
\special{pa 1028 1749}%
\special{pa 1052 1770}%
\special{pa 1078 1788}%
\special{pa 1105 1804}%
\special{pa 1134 1818}%
\special{pa 1163 1831}%
\special{pa 1194 1843}%
\special{pa 1224 1854}%
\special{pa 1240 1860}%
\special{sp}%
\put(4.4000,-14.2000){\makebox(0,0)[lb]{${\pi}_1$}}%
\put(15.6000,-9.1000){\makebox(0,0)[lb]{$y$}}%
\put(15.4000,-19.6000){\makebox(0,0)[lb]{$y-1$}}%
\put(24.9000,-14.2000){\makebox(0,0)[lb]{$x-z$}}%
\put(35.0000,-14.3000){\makebox(0,0)[lb]{$w$}}%
%
\special{pn 20}%
\special{pa 1140 990}%
\special{pa 1320 860}%
\special{fp}%
\special{sh 1}%
\special{pa 1320 860}%
\special{pa 1254 883}%
\special{pa 1277 891}%
\special{pa 1278 915}%
\special{pa 1320 860}%
\special{fp}%
%
\special{pn 20}%
\special{pa 1140 1830}%
\special{pa 1360 1900}%
\special{fp}%
\special{sh 1}%
\special{pa 1360 1900}%
\special{pa 1303 1861}%
\special{pa 1309 1884}%
\special{pa 1290 1899}%
\special{pa 1360 1900}%
\special{fp}%
\end{picture}%
\label{CPB44}
\caption{The transformations described in Theorem \ref{th:2.1} define a representation of the affine Weyl group of type $B_4^{(1)}$, that is, they satisfy the following relations: ${s_0}^2={s_1}^2=\dots={s_4}^2={{\pi}_1}^2={{\pi}_2}^2=1, \ (s_0s_1)^2=(s_0s_3)^2=(s_0s_4)^2=(s_1s_3)^2=(s_1s_4)^2=(s_2s_4)^2=1, \ (s_0s_2)^3=(s_1s_2)^3=(s_2s_3)^3=(s_3s_4)^4=1$. The symbol in each circle denotes the invariant divisor $f_i$ of the system \eqref{eq:2} (see Proposition \ref{prop:2}).}
\end{figure}
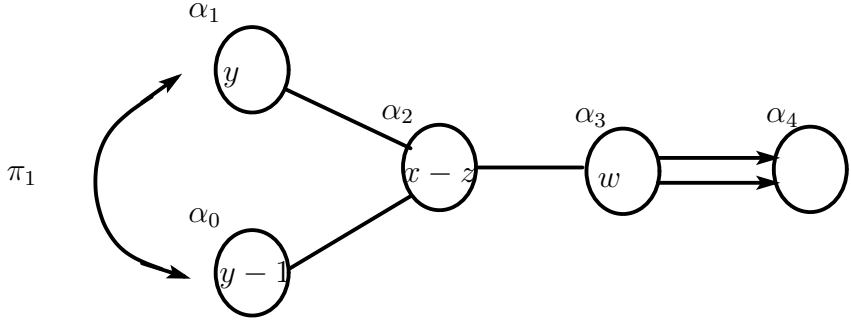

\begin{align}
\begin{split}
s_{0}: (*) \rightarrow &\left(x+\frac{\alpha_0}{y-1},y,z,w,t;-\alpha_0,\alpha_1,\alpha_2+\alpha_0,\alpha_3,\alpha_4 \right),\\
s_{1}: (*) \rightarrow &\left(x+\frac{\alpha_1}{y},y,z,w,t;\alpha_0,-\alpha_1,\alpha_2+\alpha_1,\alpha_3,\alpha_4 \right),\\
s_{2}: (*) \rightarrow &\left(x,y-\frac{\alpha_2}{x-z},z,w+\frac{\alpha_2}{x-z},t;\alpha_0+\alpha_2,\alpha_1+\alpha_2,-\alpha_2,\alpha_3+\alpha_2,\alpha_4 \right),\\
s_{3}: (*) \rightarrow &\left(x,y,z+\frac{\alpha_3}{w},w,t;\alpha_0,\alpha_1,\alpha_2+\alpha_3,-\alpha_3,\alpha_4+\alpha_3 \right),\\
s_{4}: (*) \rightarrow &\left(x,y,z,w-\frac{2\alpha_4}{z}+\frac{t}{z^2},-t;\alpha_0,\alpha_1,\alpha_2,\alpha_3+2\alpha_4,-\alpha_4 \right),\\
{\pi}_1: (*) \rightarrow &(-x,1-y,-z,-w,-t;\alpha_1,\alpha_0,\alpha_2,\alpha_3,\alpha_4),\\
{\pi}_2: (*) \rightarrow &\left(\frac{t}{z},-\frac{z}{t}(zw+\alpha_3),\frac{t}{x},-\frac{x}{t}(xy+\alpha_1),t;2\alpha_4+\alpha_3,\alpha_3,\alpha_2,\alpha_1,(\alpha_0-\alpha_1)/2 \right).
\end{split}
\end{align}
\end{Theorem}

\begin{Proposition}
Let us define the following translation operators
$$
T_1:=s_4{\pi}s_1s_2s_4s_3s_4s_3s_2s_1, \quad T_2:=s_0T_1s_0, \quad T_3:=s_2T_2s_2, \quad T_4:=s_3T_3s_3.
$$
These translation operators act on parameters $\alpha_i$ as follows$:$
\begin{align}
\begin{split}
T_1(\alpha_0,\alpha_1,\ldots,\alpha_4)=&(\alpha_0,\alpha_1,\alpha_2,\alpha_3,\alpha_4)+(1,-1,0,0,0),\\
T_2(\alpha_0,\alpha_1,\ldots,\alpha_4)=&(\alpha_0,\alpha_1,\alpha_2,\alpha_3,\alpha_4)+(-1,-1,1,0,0),\\
T_3(\alpha_0,\alpha_1,\ldots,\alpha_4)=&(\alpha_0,\alpha_1,\alpha_2,\alpha_3,\alpha_4)+(0,0,-1,1,0),\\
T_4(\alpha_0,\alpha_1,\ldots,\alpha_4)=&(\alpha_0,\alpha_1,\alpha_2,\alpha_3,\alpha_4)+(0,0,0,-1,1).
\end{split}
\end{align}
\end{Proposition}

\begin{figure}[ht]
\unitlength 0.1in
\begin{picture}(54.12,15.40)(3.00,-17.20)
%
\special{pn 20}%
\special{ar 758 407 213 227  1.5707963 6.2831853}%
\special{ar 758 407 213 227  0.0000000 1.5332553}%
%
\special{pn 20}%
\special{ar 758 1483 213 226  1.5707963 6.2831853}%
\special{ar 758 1483 213 226  0.0000000 1.5332553}%
%
\special{pn 20}%
\special{ar 2891 1498 184 222  1.5707963 6.2831853}%
\special{ar 2891 1498 184 222  0.0000000 1.5219222}%
%
\special{pn 20}%
\special{pa 946 503}%
\special{pa 1291 804}%
\special{fp}%
%
\special{pn 20}%
\special{pa 954 1399}%
\special{pa 1310 1125}%
\special{fp}%
%
\special{pn 20}%
\special{pa 1678 983}%
\special{pa 1998 983}%
\special{fp}%
%
\special{pn 20}%
\special{pa 2372 846}%
\special{pa 2732 569}%
\special{fp}%
%
\special{pn 20}%
\special{pa 2348 1127}%
\special{pa 2717 1424}%
\special{fp}%
%
\put(3.0000,-9.8100){\makebox(0,0)[lb]{}}%
\put(5.7000,-15.5000){\makebox(0,0)[lb]{$y+t$}}%
\put(12.5000,-10.4000){\makebox(0,0)[lb]{$x-z$}}%
\put(20.7500,-10.7000){\makebox(0,0)[lb]{$w$}}%
\put(27.6700,-15.9000){\makebox(0,0)[lb]{$z$}}%
\put(27.0000,-5.2000){\makebox(0,0)[lb]{$z-1$}}%
\put(6.7100,-4.9700){\makebox(0,0)[lb]{$y$}}%
%
\special{pn 20}%
\special{pa 4350 973}%
\special{pa 4629 973}%
\special{fp}%
\put(47.2300,-10.6600){\makebox(0,0)[lb]{$w$}}%
\put(39.8100,-10.3100){\makebox(0,0)[lb]{$x-z$}}%
%
\special{pn 20}%
\special{pa 2921 689}%
\special{pa 3273 689}%
\special{fp}%
\special{sh 1}%
\special{pa 3273 689}%
\special{pa 3206 669}%
\special{pa 3220 689}%
\special{pa 3206 709}%
\special{pa 3273 689}%
\special{fp}%
%
\special{pn 20}%
\special{pa 2927 1200}%
\special{pa 3281 1200}%
\special{fp}%
\special{sh 1}%
\special{pa 3281 1200}%
\special{pa 3214 1180}%
\special{pa 3228 1200}%
\special{pa 3214 1220}%
\special{pa 3281 1200}%
\special{fp}%
%
\special{pn 20}%
\special{ar 3531 1497 202 220  1.5707963 6.2831853}%
\special{ar 3531 1497 202 220  0.0000000 1.5262713}%
%
\special{pn 20}%
\special{pa 3675 551}%
\special{pa 3998 808}%
\special{fp}%
%
\special{pn 20}%
\special{pa 3708 1422}%
\special{pa 3989 1085}%
\special{fp}%
%
\special{pn 20}%
\special{pa 5024 878}%
\special{pa 5314 878}%
\special{fp}%
\special{sh 1}%
\special{pa 5314 878}%
\special{pa 5247 858}%
\special{pa 5261 878}%
\special{pa 5247 898}%
\special{pa 5314 878}%
\special{fp}%
%
\special{pn 20}%
\special{pa 5016 1067}%
\special{pa 5322 1067}%
\special{fp}%
\special{sh 1}%
\special{pa 5322 1067}%
\special{pa 5255 1047}%
\special{pa 5269 1067}%
\special{pa 5255 1087}%
\special{pa 5322 1067}%
\special{fp}%
\put(34.7900,-5.4100){\makebox(0,0)[lb]{$y$}}%
\put(33.7800,-15.6600){\makebox(0,0)[lb]{$y-1$}}%
%
\special{pn 20}%
\special{ar 1452 960 214 227  1.5707963 6.2831853}%
\special{ar 1452 960 214 227  0.0000000 1.5287650}%
%
\special{pn 20}%
\special{ar 4170 952 202 220  1.5707963 6.2831853}%
\special{ar 4170 952 202 220  0.0000000 1.5262713}%
%
\special{pn 20}%
\special{ar 4827 972 202 220  1.5707963 6.2831853}%
\special{ar 4827 972 202 220  0.0000000 1.5262713}%
%
\special{pn 20}%
\special{ar 2189 978 184 224  1.5707963 6.2831853}%
\special{ar 2189 978 184 224  0.0000000 1.5327712}%
%
\special{pn 20}%
\special{ar 3530 400 202 220  1.5707963 6.2831853}%
\special{ar 3530 400 202 220  0.0000000 1.5262713}%
%
\special{pn 20}%
\special{ar 5510 970 202 220  1.5707963 6.2831853}%
\special{ar 5510 970 202 220  0.0000000 1.5262713}%
%
\special{pn 20}%
\special{ar 2880 430 202 220  1.5707963 6.2831853}%
\special{ar 2880 430 202 220  0.0000000 1.5262713}%
\end{picture}%
\label{CPB41}
\caption{The process from the left hand side to the right hand side denotes the confluence process of the system \eqref{eq:1} to \eqref{eq:2}.}
\end{figure}
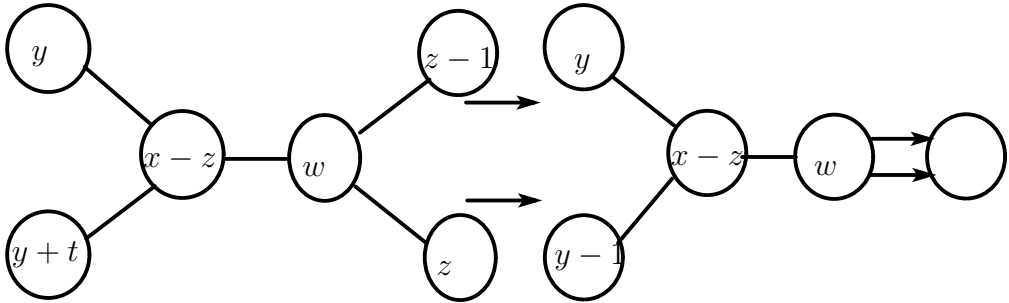

\begin{Theorem}
Let us consider a polynomial Hamiltonian system with Hamiltonian $H \in {\Bbb C}(t)[x,y,z,w]$. We assume that

$(A1)$ $deg(H)=6$ with respect to $x,y,z,w$.

$(A2)$ This system becomes again a polynomial Hamiltonian system in each coordinate system $(x_i,y_i,z_i,w_i) \ (i=0,1,\dots ,4)${\rm : \rm}
\begin{align}
\begin{split}
&x_0=1/x, \ y_0=-((y-1)x+\alpha_0)x, \ z_0=z, \ w_0=w,\\
&x_1=1/x, \ y_1=-(yx+\alpha_1)x, \ z_1=z, \ w_1=w,\\
&x_2=-((x-z)y-\alpha_2)y, \ y_2=1/y, \ z_2=z, \ w_2=w+y,\\
&x_3=x, \ y_3=y, \ z_3=1/z, \ w_3=-(wz+\alpha_3)z, \\
&x_4=x, \ y_4=y, \ z_4=z, \ w_4=w-\frac{2\alpha_4}{z}+\frac{t}{z^2}.
\end{split}
\end{align}
Then such a system coincides with the system \eqref{eq:2}.
\end{Theorem}
We remark that the B{\"a}cklund transformations $s_0,s_1,s_2,s_3$ have Noumi-Yamada's universal description for $B_4^{(1)}$ root system. But the transformation $s_4$ do not have so. In this vein, it is still an open question whether our system of type $B_4^{(1)}$ can be obtained by similarity reduction of a Drinfeld-Sokolov hierarchy.

\begin{Theorem}\label{th:2.2}
For the system \eqref{eq:1} of type $D_5^{(1)}$, we make the change of parameters and variables
\begin{gather}
\begin{gathered}\label{13}
\alpha_0=A_0, \ \alpha_1=A_1, \ \alpha_2=A_2, \ \alpha_3=A_3, \ {\alpha_4}=2A_4-\frac{1}{{\varepsilon}}, \ \alpha_5=\frac{1}{{\varepsilon}},
\end{gathered}\\
\begin{gathered}\label{14}
t=-{\varepsilon}T, \ x=1+\frac{X}{\varepsilon T}, \ y=\varepsilon TY, \ z=1+\frac{Z}{\varepsilon T}, \ w=\varepsilon TW
\end{gathered}
\end{gather}
from $\alpha_0,\alpha_1,\dots,\alpha_5,t,x,y,z,w$ to $A_0,A_1,\dots,A_4,\varepsilon,T,X,Y,Z,W$. Then the system \eqref{eq:1} can also be written in the new variables $T,X,Y,Z,W$ and parameters $A_0,A_1,\dots,\\
A_4,\varepsilon$ as a Hamiltonian system. This new system tends to the system \eqref{eq:2} of type $B_4^{(1)}$ as $\varepsilon \ra 0$.
\end{Theorem}

By proving the following theorem, we see how the degeneration process in Theorem \ref{th:2.2} works on the B{\"a}cklund transformation group $W(D_5^{(1)})=<s_0,s_1,\dots,s_5>$ described in Theorem \ref{th:1.1}.
\begin{Theorem}\label{th:2.3}
For the degeneration process in Theorem \ref{th:2.2}, we can choose a subgroup $W_{D_5^{(1)} \ra B_4^{(1)}}$ of the B{\"a}cklund transformation group $W(D_5^{(1)})$ so that $W_{D_5^{(1)} \ra B_4^{(1)}}$ converges to $W(B_4^{(1)})$ as $\varepsilon \rightarrow 0$.
\end{Theorem}

\section{ The system of type $D_4^{(2)}$ }

In this section, we present a 3-parameter family of polynomial Hamiltonian systems in dimension four explicitly given by
\begin{equation}\label{eq:4}
  \left\{
  \begin{aligned}
   \frac{dx}{dt} &=\frac{2x^2y+(1-2\beta_1-2\beta_2-2\beta_3)x+2z(zw+\beta_2)}{t}+1,\\
   \frac{dy}{dt} &=\frac{-2xy^2-(1-2\beta_1-2\beta_2-2\beta_3)y+1}{t},\\
   \frac{dz}{dt} &=\frac{2z^2w+(1-2\beta_3)z}{t}+1+\frac{2yz^2}{t},\\
   \frac{dw}{dt} &=\frac{-2zw^2-(1-2\beta_3)w-2\beta_2y-4yzw}{t}
   \end{aligned}
  \right. 
\end{equation}
with the Hamiltonian
\begin{align}
\begin{split}
H_{D_4^{(2)}}=&\frac{x^2y^2+(1-2\beta_1-2\beta_2-2\beta_3)xy-x}{t}+y\\
&+\frac{z^2w^2+(1-2\beta_3)zw}{t}+w+\frac{2yz(zw+\beta_2)}{t}.
\end{split}
\end{align}
Here $x,y,z$ and $w$ denote unknown complex variables and $\beta_1,\beta_2,\dots,\beta_4$ are complex parameters satisfying the relation:
$$
\beta_1+\beta_2+\beta_3+\beta_4=\frac{1}{2}.
$$

\begin{Proposition}\label{prop:3}
The system \eqref{eq:4} has the following invariant divisors\rm{:\rm}
\begin{center}
\begin{tabular}{|c|c|c|} \hline
codimension & invariant divisors & parameter's relation \\ \hline
1 & $f_0:=x-z$ & $\beta_1=0$  \\ \hline
1 & $f_1:=w$ & $\beta_2=0$  \\ \hline
\end{tabular}
\end{center}
\end{Proposition}

\begin{Theorem}\label{th:4.1}
The system \eqref{eq:4} admits affine Weyl group symmetry of type $D_4^{(2)}$ as the group of its B{\"a}cklund transformations, whose generators $w_1,\dots,w_4$ defined as follows$:$ with {\it the notation} $(*):=(x,y,z,w,t;\beta_1,\dots,\beta_4)$,

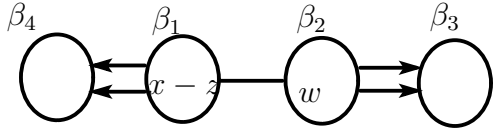
\begin{figure}[h]
\unitlength 0.1in
\begin{picture}(25.07,6.64)(2.70,-11.34)
%
\put(3.0000,-9.8900){\makebox(0,0)[lb]{}}%
%
\special{pn 20}%
\special{ar 1180 889 188 224  1.5707963 6.2831853}%
\special{ar 1180 889 188 224  0.0000000 1.5282688}%
%
\special{pn 20}%
\special{ar 1901 898 188 224  1.5707963 6.2831853}%
\special{ar 1901 898 188 224  0.0000000 1.5335795}%
%
\special{pn 20}%
\special{pa 1368 898}%
\special{pa 1697 898}%
\special{fp}%
\put(17.8000,-10.1000){\makebox(0,0)[lb]{$w$}}%
\put(10.0000,-9.8000){\makebox(0,0)[lb]{$x-z$}}%
%
\special{pn 20}%
\special{pa 984 818}%
\special{pa 726 818}%
\special{fp}%
\special{sh 1}%
\special{pa 726 818}%
\special{pa 793 838}%
\special{pa 779 818}%
\special{pa 793 798}%
\special{pa 726 818}%
\special{fp}%
%
\special{pn 20}%
\special{pa 984 955}%
\special{pa 726 955}%
\special{fp}%
\special{sh 1}%
\special{pa 726 955}%
\special{pa 793 975}%
\special{pa 779 955}%
\special{pa 793 935}%
\special{pa 726 955}%
\special{fp}%
%
\special{pn 20}%
\special{pa 2100 950}%
\special{pa 2390 950}%
\special{fp}%
\special{sh 1}%
\special{pa 2390 950}%
\special{pa 2323 930}%
\special{pa 2337 950}%
\special{pa 2323 970}%
\special{pa 2390 950}%
\special{fp}%
%
\special{pn 20}%
\special{ar 2590 910 187 224  1.5707963 6.2831853}%
\special{ar 2590 910 187 224  0.0000000 1.5333806}%
%
\special{pn 20}%
\special{pa 2100 830}%
\special{pa 2390 830}%
\special{fp}%
\special{sh 1}%
\special{pa 2390 830}%
\special{pa 2323 810}%
\special{pa 2337 830}%
\special{pa 2323 850}%
\special{pa 2390 830}%
\special{fp}%
\put(2.7000,-6.4000){\makebox(0,0)[lb]{$\beta_4$}}%
\put(10.2000,-6.6000){\makebox(0,0)[lb]{$\beta_1$}}%
\put(17.7000,-6.5000){\makebox(0,0)[lb]{$\beta_2$}}%
\put(24.6000,-6.4000){\makebox(0,0)[lb]{$\beta_3$}}%
%
\special{pn 20}%
\special{ar 530 880 187 224  1.5707963 6.2831853}%
\special{ar 530 880 187 224  0.0000000 1.5333806}%
\end{picture}%
\label{CPB42}
\caption{The transformations described in Theorem \ref{th:4.1} define a representation of the affine Weyl group of type $D_4^{(2)}$, that is, they satisfy the following relations: ${w_1}^2={w_2}^2={w_3}^2={w_4}^2=1, \ (w_1w_3)^2=(w_3w_4)^2=(w_2w_4)^2=(w_1w_2)^3=(w_1w_4)^4=(w_2w_3)^4=1$. The symbol in each circle denotes the invariant divisor $f_i$ of the system \eqref{eq:4} (see Proposition \ref{prop:3}).}
\end{figure}

\begin{align}
\begin{split}
w_{1}: (*) \rightarrow &\left(x,y-\frac{\beta_1}{x-z},z,w+\frac{\beta_1}{x-z},t;-\beta_1,\beta_2+\beta_1,\beta_3,\beta_4+\beta_1 \right),\\
w_{2}: (*) \rightarrow &\left(x,y,z+\frac{\beta_2}{w},w,t;\beta_1+\beta_2,-\beta_2,\beta_3+\beta_2,\beta_4 \right),\\
w_{3}: (*) \rightarrow &\left(x,y,z,w-\frac{2\beta_3}{z}+\frac{t}{z^2},-t;\beta_1,\beta_2+2\beta_3,-\beta_3,\beta_4 \right),\\
w_{4}: (*) \rightarrow &\left(-x-\frac{2\beta_4}{y}+\frac{1}{y^2},-y,-z,-w,-t;\beta_1+2\beta_4,\beta_2,\beta_3,-\beta_4 \right).
\end{split}
\end{align}
\end{Theorem}
We remark that the B{\"a}cklund transformations $w_1,w_2$ have Noumi-Yamada's universal description for $D_4^{(2)}$ root system. But the transformations $w_3,w_4$ do not have so. In this vein, it is still an open question whether our system of type $D_4^{(2)}$ can be obtained by similarity reduction of a Drinfeld-Sokolov hierarchy.

\begin{Proposition}
Let us define the following translation operators
\begin{align}
\begin{split}
T_1:&=w_1w_4w_2w_1w_2w_3w_1w_2,\\
T_2:&=w_3w_1w_2w_1w_4w_1w_2w_1,\\
T_3:&=w_2w_1w_4w_1w_2w_1w_3w_4w_2w_1w_2w_3w_1w_4w_1w_2w_3w_2w_1w_2.
\end{split}
\end{align}
These translation operators act on parameters $\beta_i$ as follows$:$
\begin{align}
\begin{split}
T_1(\beta_1,\beta_2,\beta_3,\beta_4)=&(\beta_1,\beta_2,\beta_3,\beta_4)+(1,-1,0,0),\\
T_2(\beta_1,\beta_2,\beta_3,\beta_4)=&(\beta_1,\beta_2,\beta_3,\beta_4)+(0,-1,1,0),\\
T_3(\beta_1,\beta_2,\beta_3,\beta_4)=&(\beta_1,\beta_2,\beta_3,\beta_4)+(0,0,-1,1).
\end{split}
\end{align}
\end{Proposition}

\begin{Theorem}
Let us consider a polynomial Hamiltonian system with Hamiltonian $H \in {\Bbb C}(t)[x,y,z,w]$. We assume that

$(A1)$ $deg(H)=6$ with respect to $x,y,z,w$.

$(A2)$ This system becomes again a polynomial Hamiltonian system in each coordinate system $(x_i,y_i,z_i,w_i) \ (i=1,2,3,4)${\rm : \rm}
\begin{align}
\begin{split}
&x_1=-((x-z)y-\beta_1)y, \ y_1=1/y, \ z_1=z, \ w_1=w+y,\\
&x_2=x, \ y_2=y, \ z_2=1/z, \ w_2=-(wz+\beta_2)z, \\
&x_3=x, \ y_3=y, \ z_3=z, \ w_3=w-\frac{2\beta_3}{z}+\frac{t}{z^2},\\
&x_4=x+\frac{2\beta_4}{y}-\frac{1}{y^2}, \ y_4=y, \ z_4=z, \ w_4=w.
\end{split}
\end{align}
Then such a system coincides with the system \eqref{eq:4}.
\end{Theorem}

\section{ The system of type $B_3^{(1)}$ }

In this section, we present a 3-parameter family of polynomial Hamiltonian systems in dimension four explicitly given by
\begin{equation}\label{eq:33333}
  \left\{
  \begin{aligned}
   \frac{dx}{dt} &=\frac{2x^2y}{t}-\frac{(\alpha_0+\alpha_1+2\alpha_2-1)x}{t}+1-\frac{2w}{t},\\
   \frac{dy}{dt} &=-\frac{2xy^2}{t}+\frac{(\alpha_0+\alpha_1+2\alpha_2-1)y}{t}+\frac{1}{t},\\
   \frac{dz}{dt} &=\frac{2z^2w}{t}-z^2+\frac{(\alpha_0+\alpha_1-1)z}{t}-\frac{2y}{t},\\
   \frac{dw}{dt} &=-\frac{2zw^2}{t}+2zw-\frac{(\alpha_0+\alpha_1-1)w}{t}+\alpha_1
   \end{aligned}
  \right. 
\end{equation}
with the Hamiltonian
\begin{align}
\begin{split}
H_{B_3^{(1)}}=&\frac{x^2y^2}{t}-\frac{(\alpha_0+\alpha_1+2\alpha_2-1)xy}{t}-\frac{x}{t}+y\\
&+\frac{z^2w^2}{t}-z^2 w+\frac{(\alpha_0+\alpha_1-1)zw}{t}-\alpha_1 z-\frac{2yw}{t}.
\end{split}
\end{align}
Here $x,y,z$ and $w$ denote unknown complex variables and $\alpha_0,\alpha_1,\alpha_2,\alpha_3$ are complex parameters satisfying the relation:
$$
\alpha_0+\alpha_1+2\alpha_2+2\alpha_3=1.
$$

\begin{Theorem}
The system \eqref{eq:33333} admits extended affine Weyl group symmetry of type $B_3^{(1)}$ as the group of its B{\"a}cklund transformations, whose generators $s_0,\dots,s_3,\pi$ defined as follows$:$ with {\it the notation} $(*):=(x,y,z,w,t;\alpha_0,\dots,\alpha_3)$,
\begin{align}
\begin{split}
s_{0}: (*) \rightarrow & \left(x,y,z+\frac{\alpha_0}{w-t},w,t;-\alpha_0,\alpha_1,\alpha_2+\alpha_0,\alpha_3 \right),\\
s_{1}: (*) \rightarrow & \left(x,y,z+\frac{\alpha_1}{w},w,t;\alpha_0,-\alpha_1,\alpha_2+\alpha_1,\alpha_3 \right),\\
s_{2}: (*) \rightarrow & \left(x,y-\frac{\alpha_2 z}{xz-1},z,w-\frac{\alpha_2 x}{xz-1},t;\alpha_0+\alpha_2,\alpha_1+\alpha_2,-\alpha_2,\alpha_3+\alpha_2 \right),\\
s_{3}: (*) \rightarrow & \left(-x-\frac{2\alpha_3}{y}+\frac{1}{y^2},-y,-z,-w,-t;\alpha_0,\alpha_1,\alpha_2+2\alpha_3,-\alpha_3 \right),\\
\pi: (*) \rightarrow &(x,y,z,w-t,-t;\alpha_1,\alpha_0,\alpha_2,\alpha_3).
\end{split}
\end{align}
\end{Theorem}

\begin{Proposition}
For the system \eqref{eq:4}, we make the change of parameters and variables
\begin{gather}
\begin{gathered}\label{PB3A3}
\alpha_0=\beta_2+2\beta_3,  \quad \alpha_1=\beta_2, \quad \alpha_2=\beta_1, \quad \alpha_3=\beta_4,
\end{gathered}\\
\begin{gathered}\label{VB3A3}
X:=x, \quad Y:=y, \quad Z:=\frac{1}{z}, \quad W:=-(zw+\beta_2)z
\end{gathered}\\
\begin{gathered}\label{Re}
W_{B_3^{(1)}}:=\{<s_0,s_1,s_2,s_3>|s_0:=w_3w_2w_3, \quad s_1:=w_2, \quad s_2:=w_1, \quad s_3:=w_4\}
\end{gathered}
\end{gather}
from $\beta_1,\beta_2,\beta_3,\beta_4,t,x,y,z,w$ to $\alpha_0,\alpha_1,\alpha_2,\alpha_3,t,X,Y,Z,W$. Then the system \eqref{eq:4} can also be written in the new variables $t,X,Y,Z,W$ and parameters $\alpha_0,\alpha_1,\alpha_2,\alpha_3$ as a Hamiltonian system. This new system tends to the Hamiltonian system \eqref{eq:33333}.
\end{Proposition}

\section{Proof of Theorems \ref{th:1.3},\ref{D5A4} and \ref{th:2.3}}

As is well-known, the degeneration from $P_{V}$ to $P_{IV}$ (see \cite{T2,T3}) is given by
$$
\alpha_0=A_0+\frac{1}{2}{\varepsilon}^{-2}, \ \alpha_1=A_1, \ \alpha_2=A_2, \ \alpha_3=-\frac{1}{2}{\varepsilon}^{-2},
$$
$$
t=\frac{1}{2}{\varepsilon}^{-2}(1+2\varepsilon T), \ x=-\frac{\varepsilon X}{1-\varepsilon X}, \ y=-\varepsilon^{-1}(1-\varepsilon X)[Y-\varepsilon (A_1+XY)].
$$
As the fourth-order analogue of the above confluence process, we consider the following coupling confluence process from the system \eqref{eq:1} by taking the above process for each coordinate system $(x,y)$ and $(z,w)$ in \eqref{eq:1}, respectively. If we take the following coupling confluence process $P_{V} \ra P_{IV}$ for each coordinate system $(x,y)$ and $(z,w)$ in \eqref{eq:1}
$$
\alpha_0=A_0-A_2-A_3+\frac{1}{2}{\varepsilon}^{-2}, \ \alpha_1=A_1, \ \alpha_2=A_2, \ \alpha_3=A_3, \ \alpha_4=-\frac{1}{2}{\varepsilon}^{-2}, \ \alpha_5=A_4,
$$
$$
t=\frac{1}{2}{\varepsilon}^{-2}(1+2\varepsilon T), \ x=-\frac{\varepsilon X}{1-\varepsilon X}, \ y=-\varepsilon^{-1}(1-\varepsilon X)[Y-\varepsilon (A_1+XY)],
$$
$$z=-\frac{\varepsilon Z}{1-\varepsilon Z}, \ w=-\varepsilon^{-1}(1-\varepsilon Z)[W-\varepsilon (A_3+ZW)],
$$
and take the limit $\varepsilon \rightarrow 0$, then we can obtain the system of type $A_4^{(1)}$ explicitly given by
\begin{equation}\label{eq:5}
  \left\{
  \begin{aligned}
   \frac{dX}{dT} &=-X^2+4XY+4ZW-2TX-2A_2-2A_4,\\
   \frac{dY}{dT} &=-2Y^2+2XY+2TY+A_1,\\
   \frac{dZ}{dT} &=-Z^2+4ZW+4YZ-2TZ-2A_4,\\
   \frac{dW}{dT} &=-2W^2+2ZW-4YW+2TW+A_3
   \end{aligned}
  \right. 
\end{equation}
with the Hamiltonian $H$
\begin{align}
\begin{split}
H_{A_4^{(1)}}=&-X^2Y+2XY^2-2TXY-(2A_2+2A_4)Y-A_1X\\
&-Z^2W+2ZW^2-2TZW-2A_4W-A_3Z+4YZW\\
&\quad (A_0+A_1+A_2+A_3+A_4=1).
\end{split}
\end{align}

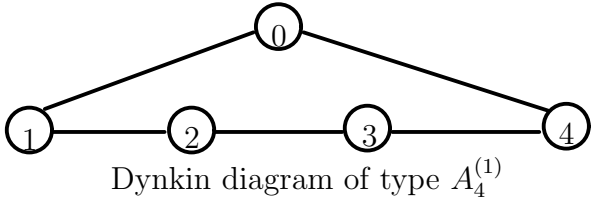
\begin{figure}[h]
\unitlength 0.1in
\begin{picture}(30.14,8.37)(19.43,-20.00)
%
\special{pn 20}%
\special{ar 2060 1820 117 117  0.0000000 6.2831853}%
%
\special{pn 20}%
\special{ar 2900 1820 117 117  0.0000000 6.2831853}%
%
\special{pn 20}%
\special{ar 3810 1810 117 117  0.0000000 6.2831853}%
%
\special{pn 20}%
\special{ar 4840 1800 117 117  0.0000000 6.2831853}%
%
\special{pn 20}%
\special{ar 3350 1280 117 117  0.0000000 6.2831853}%
%
\special{pn 20}%
\special{pa 3220 1310}%
\special{pa 2140 1710}%
\special{fp}%
%
\special{pn 20}%
\special{pa 2180 1830}%
\special{pa 2760 1830}%
\special{fp}%
%
\special{pn 20}%
\special{pa 3020 1830}%
\special{pa 3680 1830}%
\special{fp}%
%
\special{pn 20}%
\special{pa 3940 1830}%
\special{pa 4700 1830}%
\special{fp}%
%
\special{pn 20}%
\special{pa 3480 1310}%
\special{pa 4750 1720}%
\special{fp}%
\put(33.1000,-13.8000){\makebox(0,0)[lb]{0}}%
\put(20.2000,-19.2000){\makebox(0,0)[lb]{1}}%
\put(28.6000,-19.1000){\makebox(0,0)[lb]{2}}%
\put(37.8000,-19.0000){\makebox(0,0)[lb]{3}}%
\put(48.0000,-18.9000){\makebox(0,0)[lb]{4}}%
\put(24.8000,-21.7000){\makebox(0,0)[lb]{Dynkin diagram of type $A_4^{(1)}$}}%
\end{picture}%
\end{figure}

\begin{Remark}
The system \eqref{eq:5} admits extended affine Weyl group symmetry of type $A_4^{(1)}$ as the group of its B{\"a}cklund transformations, whose generators $s_0,s_1,..,s_4$ defined as follows$:$ with {\it the notation} $(*):=(X,Y,Z,W,T;A_0,A_1,\dots,A_4)$,
\begin{align}
\begin{split}
s_0: (*) \rightarrow &(X-\frac{2A_0}{X-2Y-2W+2T},Y-\frac{A_0}{X-2Y-2W+2T},\\
&Z-\frac{2A_0}{X-2Y-2W+2T},W,T;-A_0,A_1+A_0,A_2,A_3,A_4+A_0),\\
s_1: (*) \rightarrow &\left(X+\frac{A_1}{Y},Y,Z,W,T;A_0+A_1,-A_1,A_2+A_1,A_3,A_4 \right),\\
s_2: (*) \rightarrow &\left(X,Y-\frac{A_2}{X-Z},Z,W+\frac{A_2}{X-Z},T;A_0,A_1+A_2,-A_2,A_3+A_2,A_4 \right),\\
s_3: (*) \rightarrow &\left(X,Y,Z+\frac{A_3}{W},W,T;A_0,A_1,A_2+A_3,-A_3,A_4+A_3 \right),\\
s_4: (*) \rightarrow &\left(X,Y,Z,W-\frac{A_4}{Z},T;A_0+A_4,A_1,A_2,A_3+A_4,-A_4 \right).
\end{split}
\end{align}
\end{Remark}

Next, let us prove Theorem \ref{D5A4}. Notice that
$$
A_0+A_1+A_2+A_3+A_4=\alpha_0+\alpha_1+2\alpha_2+2\alpha_3+\alpha_4+\alpha_5=1
$$
and the change of variables from $(x,y,z,w)$ to $(X,Y,Z,W)$ is symplectic.

Choose $S_i \ (i=0,1,\ldots,4)$ as
$$
S_0:=s_0s_2s_3s_4s_3s_2s_0, \ S_1:=s_1, \ S_2:=s_2, \ S_3:=s_3, \ S_4:=s_5,
$$
and set $W_{D_5^{(1)} \rightarrow A_4^{(1)}}=<S_0,S_1,\ldots,S_4>$. Then we immediately have
\begin{align*}
S_0(A_0,A_1,\ldots,A_4)&=(-A_0,A_1+A_0,A_2,A_3,A_4+A_0),\\
S_1(A_0,A_1,\ldots,A_4)&=(A_0+A_1,-A_1,A_2+A_1,A_3,A_4),\\
S_2(A_0,A_1,\ldots,A_4)&=(A_0,A_1+A_2,-A_2,A_3+A_2,A_4),\\
S_3(A_0,A_1,\ldots,A_4)&=(A_0,A_1,A_2+A_3,-A_3,A_4+A_3),\\
S_4(A_0,A_1,\ldots,A_4)&=(A_0+A_4,A_1,A_2,A_3+A_4,-A_4).
\end{align*}
However, we see that $S_i(\varepsilon)$ have ambiguities of signature. For example, since
\begin{equation*}
S_3(\varepsilon^2)=\frac{\varepsilon^2}{1-2A_3\varepsilon^2},
\end{equation*}
we can choose any one of the two branches as $S_3(\varepsilon)$. Among such possibilities, we take a choice as
\begin{equation}\label{pr1:G2}
S_1(\varepsilon)=S_2(\varepsilon)=S_4(\varepsilon)=\varepsilon, \quad S_0(\varepsilon)=\varepsilon(1+2A_0 \varepsilon^2)^{-1/2}, \quad S_3(\varepsilon)=\varepsilon(1-2A_3 \varepsilon^2)^{-1/2}.
\end{equation}
Here considering in the category of formal power series, we make a convention that $(1-2A_3 \varepsilon^2)^{-1/2}$ is formal power series of $A_3 \varepsilon^2$ with constant term 1 according to
$$
(1+x)^c \sim 1+\sum_{n \geq 1} \binom{c}{n} x^n.
$$
By \eqref{11},\eqref{12},\eqref{pr1:G2} and the actions of $S_1,S_2,S_3,S_4$ on $X,Y,Z,W,T$, we can easily verify
\begin{align*}
S_1(X,Y,Z,W,T)&=\left(X+\frac{A_1}{Y},Y,Z,W,T \right),\\
S_2(X,Y,Z,W,T)&=\left(X,Y-\frac{A_2}{X-Z},Z,W+\frac{A_2}{X-Z},T \right),\\
S_3(X,Y,Z,W,T)&=\left(X,Y,Z+\frac{A_3}{W},W,(T+A_3 \varepsilon)(1-2A_3 \varepsilon^2)^{-1/2} \right),\\
S_4(X,Y,Z,W,T)&=\left(X,Y,Z,W-\frac{A_4}{Z},T \right).
\end{align*}
The forms of the actions $S_0=s_0s_2s_3s_4s_3s_2s_0$ on $X,Y,Z,W$ and $T$ are complicated, but we can see that
\begin{align*}
S_0(X,Y,Z,W,T)=&(X-\frac{2A_0}{X-2Y-2W+2T},Y-\frac{A_0}{X-2Y-2W+2T},\\
&Z-\frac{2A_0}{X-2Y-2W+2T},W,(T-A_0 \varepsilon)(1+2A_0 \varepsilon^2)^{-1/2}).
\end{align*}
The proof has thus been completed. \qed

Finally, let us prove Theorem \ref{th:2.3}. The degeneration process from the system \eqref{eq:1} to the system \eqref{eq:2} in Theorem \ref{th:2.2} is given by
$$
\alpha_0=A_0, \ \alpha_1=A_1, \ \alpha_2=A_2, \ \alpha_3=A_3, \ {\alpha_4}=2A_4-\frac{1}{{\varepsilon}}, \ \alpha_5=\frac{1}{{\varepsilon}},
$$
$$
t=-{\varepsilon}T, \ x=1+\frac{X}{\varepsilon T}, \ y=\varepsilon TY, \ z=1+\frac{Z}{\varepsilon T}, \ w=\varepsilon TW
$$
from $\alpha_0,\alpha_1,\dots,\alpha_5,t,x,y,z,w$ to $A_0,A_1,\dots,A_4,\varepsilon,T,X,Y,Z,W$. Notice that $A_0+A_1+2A_2+2A_3+2A_4=\alpha_0+\alpha_1+2\alpha_2+2\alpha_3+\alpha_4+\alpha_5=1$ and the change of variables from $(x,y,z,w)$ to $(X,Y,Z,W)$ is symplectic.
Choose $S_i$, $i=0,1,\dots,4$ as
$$
S_0:=s_0, \ S_1:=s_1, \ S_2:=s_2, \ S_3:=s_3, \ S_4:=s_4s_5=s_5s_4,
$$
which are reflections of
$$
A_0=\alpha_0, \ A_1=\alpha_1, \ A_2=\alpha_2, \ A_3=\alpha_3, \ A_4=\frac{\alpha_4+\alpha_5}{2} \ \  {\rm respectively \rm}.
$$
By using the notation $(*):=(A_0,A_1,\dots,A_4,\varepsilon)$, we can easily check
\begin{align}
\begin{split}
&S_0(*)=(-A_0,A_1,A_2+A_0,A_3,A_4,\varepsilon),\\
&S_1(*)=(A_0,-A_1,A_2+A_1,A_3,A_4,\varepsilon),\\
&S_2(*)=(A_0+A_2,A_1+A_2,-A_2,A_3+A_2,A_4,\varepsilon),\\
&S_3(*)=\left(A_0,A_1,A_2+A_3,-A_3,A_4+A_3,\frac{\varepsilon}{1+\varepsilon A_3} \right),\\
&S_4(*)=(A_0,A_1,A_2,A_3+2A_4,-A_4,-\varepsilon).
\end{split}
\end{align}
By the above relation, we will see that the group $<S_0,S_1,\dots,S_4>$ can be considered to be an affine Weyl group of the affine Lie algebra of type $B_4^{(1)}$ with respect to simple roots $A_0,A_1,\ldots,A_4$.

Now we investigate how the generators of $<S_0,S_1,\dots,S_4>$ act on $T,X,Y,Z$ and $W$. By using the notation $(**):=(X,Y,Z,W,T)$, we can verify
\begin{align}
\begin{split}
&S_0(**)=\left(X+\frac{A_0}{Y-1},Y,Z,W,T \right),\\
&S_1(**)=\left(X+\frac{A_1}{Y},Y,Z,W,T \right),\\
&S_2(**)=\left(X,Y-\frac{A_2}{X-Z},Z,W+\frac{A_2}{X-Z},T \right),\\
&S_3(**)=\left(X,Y,Z+\frac{A_3}{W},W,T(1+\varepsilon A_3) \right),\\
&S_4(**)=\left(X,Y,Z,\frac{T+\varepsilon TZW+Z^2W}{Z(\varepsilon T+Z)}-\frac{2A_4}{Z},-T \right).
\end{split}
\end{align}
The proof of Theorem \ref{th:2.3} has thus been completed.

{\it Acknowledgements.} The author would like to thank Y. Ohta, T. Suzuki, K. Takano, Y. Yamada and W. Rossman for useful discussions. In particular, Y. Yamada provided much stimulus for this research, and gave helpful advice and encouragement.

\end{document}